\documentclass[12pt]{article}
\usepackage{amsmath}
\usepackage{graphicx}
\usepackage{CJK}
\usepackage{multirow}
\usepackage{makecell}

\usepackage{subfigure}
\usepackage[justification=centering]{caption}
\input{amssym.tex}

\def\bc{\begin{center}}       \def\ec{\end{center}}
\def\ba{\begin{array}}        \def\ea{\end{array}}
\def\be{\begin{equation}}     \def\ee{\end{equation}}
\def\bea{\begin{eqnarray}}    \def\eea{\end{eqnarray}}
\def\beaa{\begin{eqnarray*}}  \def\eeaa{\end{eqnarray*}}

\def\mathbb{\Bbb}

\begin{document}
\baselineskip 18pt
\centerline{\Large\bf On the number of limit cycles for generic Lotka-Volterra system }
\vskip 0.2 true cm
\centerline{\Large\bf   and Bogdanov-Takens system under perturbations  }
\vskip 0.2 true cm
\centerline{\Large\bf  of piecewise smooth polynomials}
\vskip 0.3 true cm

\centerline{\bf  Shiyou Sui$^a$, Jihua Yang$^{b}$, Liqin Zhao$^{a,*}$}
 \centerline{ $^a$ School of Mathematical
Sciences, Beijing Normal University,} \centerline{Laboratory of Mathematics and Complex Systems, Ministry of
Education,} \centerline{Beijing 100875, The People's Republic of China}
 \centerline{ $^b$ Department of Mathematics and Computer Science,}
 \centerline{Ningxia Normal University, Guyuan 756000, The People's Republic of China}

\footnotetext[1]{
* Corresponding author.\\
E-mail:  zhaoliqin@bnu.edu.cn (L. Zhao), sui\_shiyou@163.com(S. Sui), Yangjh@mail.bnu.edu.cn(J. Yang).}

\vskip 0.2 true cm
\noindent{\bf Abstract} In this paper, we consider the bifurcation of limit cycles for generic L-V system ($\dot{x}=y+x^2-y^2\pm\frac{4}{\sqrt{3}}xy,~\dot{y}=-x+2xy$) and B-T system ($\dot{x}=y,~\dot{y}=-x+x^2$) under perturbations of piecewise smooth polynomials with degree $n$. Here the switching line is $y=0$.  By using Picard-Fuchs equations, we bound the number of zeros of first order Melnikov function which controls the number of limit cycles bifurcating from the center. It is proved that the upper bounds of the number of limit cycles for generic L-V system and B-T system are respectively $36n-65~(n\geq4),~37,57,93~(n=1,2,3)$ and $12n+6$.


\vskip 0.2 true cm
\centerline {\bf \large 1 Introduction and the main results}
\vskip 0.2 true cm
 One of the main problems in the qualitative theory of real planar continuous differential system is the determination of the numbers of limit cycles, which was proposed by Hilbert in 1990 [5,8]. In recent years, stimulated by non-smooth phenomena in the real world such as control systems [1], impact and friction mechanics [2], nonlinear oscillations [16], the theory of limit cycles for piecewise smooth differential systems has been developed. As in the smooth case, a very important issue is to determine the number of limit cycles and their distributions in the classical qualitative theory for non-smooth differential systems.

The main methods for  studying the number of limit cycles of piecewise smooth differential systems are Melnikov function method (Abelian integral) and averaging method. In [18], using first order Melnikov function, J. Yang and L. Zhao studied the bifurcation of limit cycles for the integrable non-Hamiltonian system $\dot{x}=-y(ax^2+1),~\dot{y}=x(ax^2+1),~x>0~(\dot{x}=-y(bx^2+1),~\dot{y}=x(bx^2+1),~x<0~)$ with $ab\neq0$ under perturbations of piecewise smooth polynomials with degree $n$ and switching line $x=0$. By the averaging theory of first order for discontinuous differential systems, in [13], J. Llibre and A. C. Mereu considered the quadratic isochronous centers $S_1:~\dot{x}=-y+x^2,~\dot{y}=x+xy$ and $S_2:~\dot{x}=-y+x^2-y^2,~\dot{y}=x+2xy$ when they are perturbed inside the class of all discontinuous quadratic polynomials with the straight line of discontinuity $y=0$, and they obtain that there are at least 4 and 5 limit cycles for $S_1$ and $S_2$, respectively. Recently, by using Picard-Fuchs equations and Chebyshev criterion, J. Yang and L. Zhao [17] got the exactly 5 and 6  limit cycles for $S_1$ and $S_2$, respectively. For more, one is recommended to see [4,9-12].

It is konwn that the systems that have quadratic centers at $(0,0)$ having all their orbits is cubic can be classified as the following ($z=x+iy$) [6,19]:\\
(i) The Hamiltonian system $Q_3^H$:
$$\dot{z}=-iz-z^2+2\left|z\right|^2+(b+ic)\bar{z}^2,$$
$$H(x,y)=\frac{1}{2}(x^2+y^2)+\left(\frac{b}{3}-1\right)x^3+cx^2y-(1+b)xy^2-\frac{c}{3}y^3.$$
(ii) The Hamiltonian triangle:
$$\dot{z}=-iz+\bar{z}^2,$$
$$H(x,y)=(1-2x)\left[\frac{1}{2}y^2-\frac{1}{6}(x+1)^2\right].$$
(iii) The reversible system:
 $$\dot{z}=-iz+(2b+1)z^2+2\left|z\right|^2+b\bar{z}^2,~~b\neq1,$$
 $$H(x,y)=X^{-3}\left[\frac{1}{2}y^2+\frac{1}{8(b+1)^2}\left(\frac{1-3b}{b+1}X^2+2\frac{b-1}{b+1}X+\frac{3-b}{3b+3}\right)\right],$$
where $X=1+2(b+1)x$ and the integrating factor $R(x,y)=X^{-4}$.\\
(iv) The generic Lotka-Volterra system:
$$\dot{z}=-iz+(1-ci)z^2+ci\overline{z}^2,~~c=\pm\frac{1}{\sqrt{3}},$$
$$H(x,y)=-\frac{1}{2}(1+2x)(1+x\pm\sqrt{3}y)\left(1+x\pm\frac{\sqrt{3}}{3}y\right)^{-3},$$
where the integrating factor $R(x,y)=\left(1+x\pm\frac{\sqrt{3}}{3}y\right)^{-4}$.

In this paper, we will consider the bifurcation of limit cycle for the generic Lotka-Volterra system
$$\left\{\begin{array}{lc}
\dot{x}=y+x^2-y^2\pm\displaystyle\frac{4}{\sqrt{3}}xy,\\
\dot{y}=-x+2xy,
\end{array}
\right.\eqno{(1.1)_{\pm}}$$
and  the Bogdanov-Takens system
$$\left\{\begin{array}{lc}
\dot{x}=y,\\
\dot{y}=-x+x^2
\end{array}
\right.\eqno{(1.2)}$$
which belongs to $Q_3^H$, when they are perturbed inside the piecewise smooth polynomials with degree n and switching line $y=0$. Using linear transformation [19], the system  $(1.1)_{\pm}$ can be reduced  to
$$\left\{\begin{array}{lc}
\dot{x}=xy,\\
\dot{y}=\displaystyle\frac{3}{2}y^2-\displaystyle\frac{9}{8}x^2+\displaystyle\frac{3}{2}x-\displaystyle\frac{3}{8}.
\end{array}
\right.\eqno{(1.3)}$$
Taking $X=2x-1,~Y=2\sqrt{2}y,~T=\frac{\sqrt{2}}{2}t$ in (1.2), we have (here we rewrite $X,Y$ as $x,y$)
 $$\left\{\begin{array}{lc}
\dot{x}=y,\\
\dot{y}=-1+x^2.
\end{array}
\right.\eqno{(1.4)}$$

The main results are the follows.

\noindent
{\bf Theorem 1.1.}~ Let $0<|\varepsilon|\ll 1$, and $f^\pm (x,y)$ and $g^\pm (x,y)$
are any polynomials of degree $n$. Consider the following  perturbations of systems (1.3) and (1.4):
$$
\left(
  \begin{array}{c}
          \dot{x} \\
          \dot{y}
 \end{array}
 \right)=\begin{cases}
 \left(
  \begin{array}{c}
       xy+\varepsilon f^+(x,y) \\
       \frac{3}{2}y^2-\frac{9}{8}x^2+\frac{3}{4}x-\frac{3}{8}+\varepsilon g^+(x,y)
 \end{array}
 \right), \quad y>0,\\
 \left(
  \begin{array}{c}
         xy+\varepsilon f^-(x,y) \\
        \frac{3}{2}y^2-\frac{9}{8}x^2+\frac{3}{4}x-\frac{3}{8}+\varepsilon g^-(x,y)
 \end{array}
 \right),\quad y<0,\\
 \end{cases}\eqno{(1.5)}
$$
$$
\left(
  \begin{array}{c}
          \dot{x} \\
          \dot{y}
 \end{array}
 \right)=\begin{cases}
 \left(
  \begin{array}{c}
          y+\varepsilon f^+(x,y) \\
          -1+x^2+\varepsilon g^+(x,y)
 \end{array}
 \right), \quad y>0,\\
 \left(
  \begin{array}{c}
         y+\varepsilon f^-(x,y) \\
          -1+x^2+\varepsilon g^-(x,y)
 \end{array}
 \right),\quad y<0.\\
 \end{cases}\eqno{(1.6)}
$$
Then, by using the first order Melnikov function in $\varepsilon$, the upper bounds (counting the
multiplicity) of number of limit cycles of systems (1.5) and (1.6) bifurcating from the
period annuluses are  $36n-65(n\ge 4)$, $37,57,93(n=1,2,3)$ and $12n+6$ respectively.

\noindent{\bf Remark 1.1.}~ Y. Zhao {\it et al} [19] and B. Li {\it et al} [7,14,15] considered respectively systems (1.3) and (1.4) under continuous perturbations of arbitrary polynomials with degree $n$. It is proved that the isolated zeros of Abelian integrals for perturbed system (1.3), taking into their multiplicities, does not exceed 7n. For perturbed system (1.4), the exactly upper bound of the first order Melnikov function (Abelian integral) is $n-1$, and the exactly upper bound of the second order Melnikov function is $2n-2$ (n is even) or $2n-3$ (n is odd) when the first order Melnikov function vanishes.

The paper is organized as follow. In section 2, we will give some preliminaries. In section 3, we will prove Theorem 1.1 for system (1.5). First we will obtain the algebraic structure of the first order Melnikov function $M(h)$. And then we prove that there exists a second-order differential operator that can simplify $M(h)$. Finally, the main results is proved by using the Chebyshev space. In section 4, Theorem 1.1 for system (1.6) is proved. The method is similar to the proof of system (1.5). The main difference is that $M(h)$ has five basic integrals for (1.5) while four for (1.6). Hence, we use some different techniques in the proof.

\vskip 0.2 true cm
\centerline {\bf \large 2 Preliminaries}
\vskip 0.2 true cm
We first introduce the first order Melnikov function of discontinuous differential systems. Consider the following system:
$$
(\dot{x},\ \dot{y})=\begin{cases}
(P^+(x,y)+\varepsilon f^+(x,y),\ Q^+(x,y)+\varepsilon g^+(x,y)),\ \ y>0,\\
(P^-(x,y)+\varepsilon g^-(x,y),\ Q^-(x,y)+\varepsilon g^-(x,y)),\ \ y<0,\\
\end{cases}\eqno(2.1)$$
where $0<|\varepsilon|\ll1$, and $f^\pm(x,y)$ and $g^\pm(x,y)$ are polynomials with degree $n$.  System (2.1) has two subsystems:
$$\left\{{\begin{split}
\dot{x}&=P^+(x,y)+\varepsilon f^+(x,y), \\
\dot{y}&=Q^+(x,y)+\varepsilon g^+(x,y),
\end{split}}~~~~~y>0,
\right.\eqno{(2.2)}$$
and
$$\left\{{\begin{split}
\dot{x}&=P^-(x,y)+\varepsilon f^-(x,y), \\
\dot{y}&=Q^-(x,y)+\varepsilon g^-(x,y),
\end{split}}~~~~~y<0.
\right.\eqno{(2.3)}$$
Suppose that  system (2.2)$_{\varepsilon=0}$ is integrable with the first integral $H^+(x,y)$ and integrating factor $\mu_1(x,y)$, and system (2.3)$_{\varepsilon=0}$ is integrable with the first integral $H^-(x,y)$ and  integrating factor $\mu_2(x,y)$.
We also suppose that (2.1)$_{\varepsilon=0}$ has a family of periodic orbits around the origin and satisfies the following two assumptions.

{\bf Assumption (I).} There exist an interval $\Sigma=(\alpha,\beta)$, and two points $A(h)=(a(h),0)$ and $B(h)=(b(h),0)$ such that for $h\in \Sigma$
$$H^+(A(h))=H^+(B(h))=h,\ \ \ H^-(A(h))=H^-(B(h))=\tilde{h},\ \ \ a(h)\neq b(h).$$

{\bf Assumption (II).} The subsystem (2.2)$_{\varepsilon=0}$ has an orbital arc $L^+_h$ starting from $A(h)$ and ending at $B(h)$ defined by $H^+(x,y)=h,\ y\geq0$.  The subsystem (2.3)$_{\varepsilon=0}$ has an orbital arc $L^-_h$ starting from $B(h)$ and ending at $A(h)$ defined by $H^-(x,y)=H^-(B(h)),\ y<0$.

 Under the Assumptions (I) and (II), (2.1)$_{\varepsilon=0}$ has a family of non-smooth periodic orbits $L_h=L^+_h\cup L^-_h(h\in \Sigma)$. For definiteness, we assume that the orbits $L_h$ for $h\in \Sigma$ orientate clockwise(see Fig.\,1).
 The authors [12] established a bifurcation function $F(h,\varepsilon)$ for (2.1).
 Let $F(h,0)=M(h)$.  In [4] and [9], the authors obtained the following results.

\noindent
{\bf Lemma 2.1}.\,[4,9] Under the assumptions (I) and (II), we have

(i)~~If $M(h)$ has $k$ zeros in $h$ on the interval $\Sigma$ with each having an odd multiplicity, then (2.2) has at least $k$ limit cycles bifurcating from the period annulus for $0<|\varepsilon|\ll 1$;

(ii)~~If $M(h)$ has at most $k$ zeros in $h$ on the interval $\Sigma$, taking into
account the multiplicity , then there exist at most $k$ limit cycles of (2.1) bifurcating from the period annulus;

(iii)~~The first order Melnikov function $M(h)$ of system (2.1) has the following form
\begin{align}
M(h)&=\frac{H_x^+(A)}{H_x^-(A)}\left[\frac{H_x^-(B)}{H_x^+(B)}\int_{L_h^+}\mu_1(x,y)g^+(x,y)
{\rm d}x-\mu_1(x,y)f^+(x,y){\rm d}y\right.\notag\\
&\left.+\int_{L_h^-}\mu_2(x,y)g^-(x,y){\rm d}x-\mu_2(x,y)f^-(x,y){\rm d}y\right].\tag{2.4}
\end{align}

Let $\mathbb{V}$ be a finite-dimensional vector space of functions, real-analytic on an open interval $\mathbb{I}$. Next, we give the relation of the number of zeros about second order linear  homogeneous equation and non-homogeneous equation which will be used in the proof of Theorem 1.1.

\noindent{\bf Definition 2.1} ([3]). We say that $\mathbb{V}$ is a {\bf Chebyshev space}, provided that each non-zero function in $\mathbb{V}$ has at most ${\rm dim}(\mathbb{V})-1$ zeros, counted with multiplicity.

Let $\mathbb{S}$ be the solution space of a second order linear analytic differential equation
$$x''+a_1(t)x'+a_2(t)x=0\eqno{(2.5)}$$
on an open interval $\mathbb{I}$.

\noindent{\bf Proposition 2.1} ([3]). The solution space $\mathbb{S}$ of (2.5) is a Chebyshev space of the interval $\mathbb{I}$ if and only if there exists a nowhere vanishing solution $x_0(t)\in \mathbb{S}~(x_0(t)\neq0,~\forall t\in \mathbb{I})$.

\noindent{\bf Proposition 2.2} ([3]). Suppose the solution space of the homogeneous equation (2.5) is a Chebyshev space and let $R(t)$ be an analytic function on $\mathbb{I}$ having $l$ zeros (counted with multiplicity). Then every solution $x(t)$ of the non-homogeneous equation
$$x''+a_1(t)x'+a_2(t)x=R(t)$$
has at most $l+2$ zeros on $\mathbb{I}$.

\vskip 0.2 true cm
\centerline {\bf \large 3 Proof of Theorem 1.1:  for system (1.5)}
\vskip 0.2 true cm

In this section, we will prove Theorem 1.1 for system (1.5). For $\varepsilon=0$, system (1.5) has the first integral
$$H(x,y)=x^{-3}\left(\frac{1}{2}y^2-\frac{9}{8}x^2+\frac{3}{4}x-\frac{1}{8}\right)
=h,~~~h\in(-\frac{1}{2},0)\eqno{(3.1)}$$
with integrating factor $\mu(x,y)=x^{-4}$. The point (1,0) is an elementary center corresponding to $h=-\frac{1}{2}$ and the point $(\frac{1}{3},0)$ is saddle corresponding to $h=0$ (see Fig.\,2).
Suppose that
$$f^{\pm}(x,y)=\sum\limits_{i+j=0}^na_{i,j}^{\pm}x^iy^j,~~g^{\pm}(x,y)=\sum\limits_{i+j=0}^nb_{i,j}^{\pm}x^iy^j.$$
By (2.4), for $h\in(-\frac{1}{2},0)$, we have
$$
M(h)=\int_{\Gamma_h^+}x^{-4}[g^+{\rm d}x-f^+{\rm d}y]
+\int_{\Gamma_h^-}x^{-4}[g^-{\rm d}x-f^-{\rm d}y],\eqno(3.2)$$
where
$$\Gamma_h^+(\Gamma_h^-)=\bigl\{\left.(x,y)\right|H(x,y)=h,y>0(y<0),
h\in(-\frac{1}{2},0)\bigr\}.$$
For $i\in\mathbb{Z}$ and $j\geq0$, denote
$$I_{i,j}(h)=\int_{\Gamma_h^+}x^{i-4}y^j{\rm d}x,~~J_{i,j}(h)=\int_{\Gamma_h^-}x^{i-4}y^j{\rm d}x.\eqno{(3.3)}$$

\noindent{\bf Lemma 3.1.}  If $n\geq4$, then for $h\in(-\frac{1}{2},0)$  $M(h)$ can be expressed as
$$M(h)=\frac{1}{h^{n-3}}\left[{\bf \sigma(h)}{\bf V_1(h)}+{\bf \tau(h)}{\bf V_2(h)}\right],$$
where
$$~~~{\bf \sigma(h)}=(\alpha_1(h), \alpha_2(h)), ~~{\bf \tau(h)}=(\alpha_3(h), \alpha_4(h),
\alpha_5(h)), \eqno(3.4)$$
$${\bf V_1(h)}=({\bf I_{0,1}}(h),{\bf I_{-1,1}}(h) )^T, ~~
{\bf V_2(h)}=({\bf I_{1,0}}(h),{\bf I_{0,0}}(h), {\bf I_{0,2}}(h) )^T,\eqno(3.5)$$
and $\alpha_i(h)$ are polynomials of $h$ with
${\rm deg}\alpha_i(h)\leq n-3$($i=1,\cdots, 5$).
Particularly,
$$M(h)=\frac{1}{h^{n-2}}\left[{\bf \sigma(h)}{\bf V_1(h)}+{\bf \tau(h)}{\bf V_2(h)}\right],$$
with ${\rm deg}\alpha_i(h)=n-2$($i=1,\cdots, 5$) for $n=3,2$;  and
$$M(h)={\bf \sigma(h)}{\bf V_1(h)}+(\alpha_3(h),\alpha_4(h))({\bf I_{1,0}}(h),{\bf I_{0,0}}(h))^T$$
with ${\rm deg}\alpha_i(h)=0$($i=1,\cdots, 4$) for $n=1$.

\vskip 0.2 true cm

\noindent{\bf Proof.} Suppose that $n\ge 4$ and
the orbit $\Gamma_h^+$ intersects the $x$-axis at points
$A(x_a(h),0)$ and $B(x_b(h),0)$ with $x_a(h)<1<x_b(h)$(see Fig. 2).
Let $D$ be the interior of $\Gamma_h^+\cup\overrightarrow{BA}$. Then, for $i\in\mathbb{Z}$ and  $j\geq 0$, applying the Green's formula, we can obtain
$$ \int_{\Gamma_h^+}x^{i-4}y^j{\rm d}y
=\oint_{\Gamma_h^+\cup\overrightarrow{BA}}x^{i-4}y^j{\rm d}y-\int_{\overrightarrow{BA}}x^{i-4}y^j{\rm d}y
=-(i-4)\iint\limits_{D}x^{i-5}y^j{\rm d}x{\rm d}y,$$
and
$$\int_{\Gamma_h^+}x^{i-5}y^{j+1}{\rm d}x=\oint_{\Gamma_h^+\cup\overrightarrow{BA}}x^{i-5}y^{j+1}{\rm d}x=(j+1)\iint\limits_{D}x^{i-5}y^j{\rm d}x{\rm d}y.$$
Hence, we get
$$\int_{\Gamma_h^+}x^{i-4}y^j{\rm d}y=-\frac{i-4}{j+1}\int_{\Gamma_h^+}x^{i-5}y^{j+1}{\rm d}x,
~~i\in\mathbb{Z}, j\geq 0. $$
Similarly, we have
$$\int_{\Gamma_h^-}x^{i-4}y^j{\rm d}y=-\frac{i-4}{j+1}\int_{\Gamma_h^-}x^{i-5}y^{j+1}{\rm d}x,~~i\in\mathbb{Z}, j\geq 0.$$
By direct computations, we have
$$J_{i,j}(h)=(-1)^{j+1}I_{i,j}(h).\eqno{(3.6)}$$
Hence, we get
\begin{align}
M(h)
&=\sum\limits_{i+j=0}^nb_{i,j}^
+I_{i,j}(h)+\sum\limits_{i+j=0}^na_{i,j}^+\frac{i-4}{j+1}I_{i-1,j+1}(h)\notag\\
\quad\quad
&+\sum\limits_{i+j=0}^nb_{i,j}^-J_{i,j}(h)+\sum\limits_{i+j=0}^n
a_{i,j}^-\frac{i-4}{j+1}J_{i-1,j+1}(h)\notag\\
&=\sum\limits_{\begin{subarray}{c}i+j=0\\i\geq-1,j\geq0\end{subarray}}^n\tilde{a}
_{i,j}I_{i,j}(h)+
\sum\limits_{\begin{subarray}{c}i+j=0\\i\geq-1,j\geq0\end{subarray}}^n\tilde{b}_
{i,j}J_{i,j}(h)
=:\sum\limits_{\begin{subarray}{c}i+j=0
\\i\geq-1,j\geq0\end{subarray}}^n\rho_{i,j}I_{i,j}(h),\notag
\end{align}
where $\tilde{a}_{i,j},\tilde{b}_{i,j},\rho_{i,j}$ are arbitrary real constants.
From (3.1), we have
$$x^{-3}y\frac{\partial{y}}{\partial{x}}-\frac{3}{2}x^{-4}y^2+\frac{9}{8}x^{-2}
-\frac{3}{2}x^{-3}+\frac{3}{8}x^{-4}=0.\eqno{(3.7)}$$
Multiplying (3.7) by $x^iy^{j-2}{\rm d}x$, integrating over $\Gamma_h^+$, we have
$$(2i+3j-6)I_{i,j}(h)=-2j[-\frac{9}{8}I_{i+2,j-2}(h)+\frac{3}{2}I_{i+1,j-2}(h)
-\frac{3}{8}I_{i,j-2}(h)],~~j\geq2.\eqno{(3.8)}$$
Hence, we get
$$M(h)=\sum\limits_{i=-1}^{n-1}c_{i,1}{\bf I_{i,1}}(h)+\sum\limits_{j=0}^nd_{j,1}{\bf I_{j,0}}(h)+\rho_{0,2}{\bf I_{0,2}}(h),\eqno{(3.9)}$$
where $c_{i,1},d_{j,1}$ are some real constants.

\vskip 0.2 true cm
(1)~ We first assert that
\begin{eqnarray*}\qquad \qquad\qquad\qquad\qquad
\begin{cases}
{\bf I_{1,1}}(h)={\bf I_{0,1}}(h),\\
{\bf I_{2,1}}(h)=\frac{1}{h}\left(-\frac{3}{2}{\bf I_{0,1}}(h)+{\bf I_{-1,1}}(h)\right)\\
{\bf I_{3,1}}(h)=\frac{1}{h}\left(-\frac{1}{2}{\bf I_{0,1}}(h)\right),\\
{\bf I_{2,0}}(h)=\frac{4}{3}{\bf I_{1,0}}(h)-\frac{1}{3}{\bf I_{0,0}}(h),\\
{\bf I_{3,0}}(h)=\frac{1}{h}\left(\frac{1}{2}{\bf I_{0,2}}(h)-\frac{3}{4}{\bf I_{1,0}}(h)+\frac{1}{4}{\bf I_{0,0}}(h)\right).
\end{cases}~~~~~~\qquad(3.10)\end{eqnarray*}
In fact, by Proposition 2.3 of [18], we have
$$\oint_{\Gamma_h^+\cup\Gamma_h^-}x^{-3}y{\rm d}x=\oint_{\Gamma_h^+\cup\Gamma_h^-}x^{-4}y{\rm d}x.$$
Since $H(x,y)=x^{-3}\left(\frac{1}{2}y^2-\frac{9}{8}x^2+\frac{3}{4}x-\frac{1}{8}\right)
=h$ is symmetrical about x-axis, it holds that
$$\int_{\Gamma_h^+}x^{-3}y{\rm d}x=\int_{\Gamma_h^-}x^{-3}y{\rm d}x,~~\int_{\Gamma_h^+}x^{-4}y{\rm d}x=\int_{\Gamma_h^-}x^{-4}y{\rm d}x.$$
Hence, ${\bf I_{1,1}}(h)={\bf I_{0,1}}(h)$.
Rewrite (3.1) in the form
$$\frac{1}{2}y^2-\frac{9}{8}x^2+\frac{3}{4}x-\frac{1}{8}=hx^3,\eqno{(3.11)}$$
which yields
$$hI_{i,j}(h)=\frac{1}{2}I_{i-3,j+2}(h)-\frac{9}{8}I_{i-1,j}(h)+\frac{3}{4}
I_{i-2,j}(h)-\frac{1}{8}I_{i-3,j}(h).\eqno{(3.12)}$$
By (3.8), we have
$$(2i+3j-6)I_{i-3,j+2}(h)=-2(j+2)[-\frac{9}{8}I_{i-1,j}(h)+\frac{3}{2}I_{i-2,j}
(h)-\frac{3}{8}I_{i-3,j}(h)].\eqno{(3.13)}$$
Substituting (3.13) into (3.12), we obtain
$$(2i+3j-6)hI_{i,j}(h)=-(i+j-4)\frac{9}{4}I_{i-1,j}(h)+(2i+j-10)\frac{3}{4}
I_{i-2,j}(h)-(i-6)\frac{1}{4}I_{i-3,j}(h).\eqno{(3.14)}$$
By direct computations, taking  $(i,j)=(2,1)$ and $(3,1)$ in (3.14), and noting that ${\bf I_{1,1}}(h)={\bf I_{0,1}}(h)$,
we can obtain the second and the third formulas in (3.10).
And taking $(i,j)=(0,2)$ in (3.8) and $(i,j)=(3,0)$ in (3.12), we have the last two formulas in (3.10), respectively.

\vskip 0.2 true cm
(2)~ Next we will prove that, for $i,j\geq4$, we have
\begin{align}
{\bf {\bf I_{i,1}}}(h)&=\frac{1}{h^{i-2}}\left[\alpha_{i,1}(h){\bf I_{0,1}}(h)+\beta_{i,1}(h){\bf I_{-1,1}}(h)\right],\tag{3.15}\\
{\bf I_{j,0}}(h)&=\frac{1}{h^{j-3}}\left[\gamma_{j,1}(h){\bf I_{1,0}}(h)+\delta_{j,1}(h){\bf I_{0,0}}(h)\right],\tag{3.16}
\end{align}
where $\alpha_{i,1}(h),\beta_{i,1}(h),\gamma_{j,1}(h),\delta_{j,1}(h)$ are polynomials of $h$ with ${\rm deg}\alpha_{i,1}(h)\leq[(2j-5)/3]$, ${\rm deg}\beta_{i,1}(h)$, ${\rm deg}\gamma_{j,1}(h)$ and ${\rm deg}\delta_{j,1}(h)\leq[(2j-7)/3]$, $[s]$ denotes the integer part of $s$.

In fact, using (3.10) and (3.14), it is easy to check that (3.15) holds for $i=4,5,6$. Suppose that (3.15) holds for $k\leq i-1$. Then, for $k=i$, it follows form (3.14) that
$${\bf {\bf I_{i,1}}}(h)=\frac{1}{h}\left[-\frac{9(i-3)}{4(2i-3)}{\bf I_{i-1,1}}(h)+\frac{3(2i-9)}{4(2i-3)}{\bf I_{i-2,1}}(h)-\frac{i-6}{4(2i-3)}{\bf I_{i-3,1}}(h)\right].$$
By induction assumption, we have
$${\begin{split}
{\bf {\bf I_{i,1}}}(h)&=\frac{1}{h^{i-2}}\left[\left(-\frac{9(i-3)}{4(2i-3)}\alpha_{i-1,1}(h)+\frac{3(2i-9)}{4(2i-3)}h\alpha_{i-2,1}(h)
-\frac{i-6}{4(2i-3)}h^2\alpha_{i-3}(h)\right){\bf I_{0,1}}(h)\right.\\
&\quad+\left.\left(-\frac{9(i-3)}{4(2i-3)}\beta_{i-1,1}(h)+\frac{3(2i-9)}{4(2i-3)}h\beta_{i-2,1}(h)
-\frac{i-6}{4(2i-3)}h^2\beta_{i-3}(h)\right){\bf I_{-1,1}}(h)\right]\\
&:=\frac{1}{h^{i-2}}\left[\alpha_{i,1}(h){\bf I_{0,1}}(h)+\beta_{i,1}(h){\bf I_{-1,1}}(h)\right],
\end{split}}$$
where
$${\rm deg}\alpha_{i,1}(h)\leq\max\left\{\left[\frac{2i-2-5}{3}\right],\left[\frac{2i-4-5}{3}\right]+1,\left[\frac{2i-6-5}{3}\right]+2\right\}=\left[\frac{2i-5}{3}\right],$$
$${\rm deg}\beta_{i,1}(h)\leq\max\left\{\left[\frac{2i-2-7}{3}\right],\left[\frac{2i-4-7}{3}\right]+1,\left[\frac{2i-6-7}{3}\right]+2\right\}=\left[\frac{2i-7}{3}\right].$$
This ends the proof of (3.15). By similar argument, we can get (3.16). This ends the proof. $\diamondsuit$

\vskip 0.2 true cm

\noindent{\bf Lemma 3.2.}  We have the following results:

(1)~The vector functions ${\bf V_1(h)}
$ and ${\bf V_2(h)}$ satisfy respectively the following Picard-Fuchs equations:
$${\bf V_1(h)}=(A_1h+B_1){\bf V_1'(h)},~~
{\bf V_2(h)}=(A_2h+B_2){\bf V_2'(h)},\eqno{(3.17)}$$
where
\begin{eqnarray*}
A_1=\left(\begin{matrix}
                \frac{3}{2}&0\\
                 \frac{3}{2}&\frac{3}{4}
                \end{matrix}\right),\ \
            B_1=\left(\begin{matrix}
                \frac{9}{8}&-\frac{3}{8}\\
                 \frac{27}{16}&-\frac{9}{16}
                \end{matrix}\right),\\
\end{eqnarray*}
\begin{eqnarray*}
A_2=\left(\begin{matrix}
                \frac{11}{2}&-1&0\\
                 10&-1&0\\
                 \frac{13}{4}&-1&1
                \end{matrix}\right),\ \
            B_2=\left(\begin{matrix}
                \frac{27}{8}&-\frac{9}{8}&0\\
                 \frac{27}{4}&-\frac{9}{4}&0\\
                 \frac{27}{16}&-\frac{9}{16}&0
                \end{matrix}\right).\\\end{eqnarray*}

(2)~ The functions ${\bf V_2(h)}$  satisfies the equation
$$
h(2h+1){\bf V_2''(h)}=
\left(\begin{matrix}
                -\frac{1}{18}(116h+171)&\frac{4}{9}h\\
                  -\frac{1}{18}(800h+513)&\frac{4}{9}h\\
                   -15(2h+1)&\frac{1}{2}(2h+1)
                \end{matrix}\right)
\left(\begin{matrix}
              {\bf I_{1,0}}'(h)\\
              {\bf I_{0,0}}'(h)
                \end{matrix}\right).\eqno{(3.18)}$$

\vskip 0.2 true cm

\noindent{\bf Proof.}
By direct calculation, we have
$$I_{i,j}'(h)=\int_{x_{a}(h)}^{x_b(h)}jx^{i-4}y^{j-1}\frac{\partial{y}}{\partial{h}}{\rm d}x+x_b(h)^{i-4}y(x_{b}(h),h)^j\frac{\partial{x_{b}(h)}}{\partial{h}}-x_a(h)^{i-4}y(x_{a}(h),h)^j\frac{\partial{x_a(h)}}{\partial{h}},$$
where $\Gamma_h^+$ intersects the x-axis at points  $(x_a(h),0)$ and $(x_b(h),0)$.
From
$$-\frac{9}{8}x_{a}(h)^2+\frac{3}{4}x_a(h)-\frac{1}{8}=hx_a(h)^3,$$
we get
$$\frac{\partial{x_a(h)}}{\partial{h}}=\frac{8x_a(h)^4}{3(3x_a(h)-1)(x_a(h)-1)}.$$
Since $x_a(h)\in(\frac{1}{3},1)$, we have $\frac{\partial{x_a(h)}}{\partial{h}}\neq\infty$. Similarly, $\frac{\partial{x_b(h)}}{\partial{h}}\neq\infty$.  Hence, by $y(x_a(h),h)=y(x_b(h),h)=0$, we can obtain that
$$I_{i,j}'(h)=\int_{x_{a}(h)}^{x_b(h)}jx^{i-4}y^{j-1}\frac{\partial{y}}{\partial{h}}{\rm d}x.$$
It follows from (3.11) that
$$\frac{\partial y}{\partial h}=\frac{x^3}{y}, $$
which implies
$$I_{i,j}'(h)=jI_{i+3,j-2}(h).$$
Therefore,
$$I_{i,j}(h)=\frac{1}{j+2}I_{i-3,j+2}'(h).\eqno{(3.19)}$$
Taking $(i,j)=(1,1),(0,1)$ in (3.14), and noting that ${\bf I_{1,1}}(h)={\bf I_{0,1}}(h)$, we have
$${\bf I_{-2,1}}(h)=\left(-\frac{4}{5}h-\frac{18}{5}\right){\bf I_{0,1}}(h)+\frac{21}{5}{\bf I_{-1,1}}(h),$$
$${\bf I_{-3,1}}(h)=\left(-\frac{28}{5}h-\frac{81}{5}\right){\bf I_{0,1}}(h)+\frac{72}{5}{\bf I_{-1,1}}(h),$$
respectively. Hence, taking $(i,j)=(-3,3)$ in (3.8), we get
$${\begin{split}
{\bf I_{-3,3}}(h)&=-\frac{9}{4}{\bf I_{-1,1}}(h)+3{\bf I_{-2,1}}(h)-\frac{3}{4}{\bf I_{-3,1}}(h)\\
&=\left(\frac{9}{5}h+\frac{27}{20}\right){\bf I_{0,1}}(h)-\frac{9}{20}{\bf I_{-1,1}}(h).
\end{split}}$$
Therefore, taking (i,j)=(0,1) in (3.19), we have
$${\bf I_{0,1}}(h)=\left(\frac{3}{2}h+\frac{9}{8}\right){\bf I_{0,1}}'(h)-\frac{3}{8}{\bf I_{-1,1}}'(h),$$
which implies the first equation of ${\bf V_1(h)}=(A_1h+B_1){\bf V_1'(h)}$. Similarly, we can obtain other equations of (3.17).

 Differentiating both side of the second equation of (3.17) yields
$$(A_2h+B_2){\bf V_2''(h)}=(E-A_2){\bf V_2'(h)},$$
where $E$ is a $3\times3$ identity matrix, which means (3.18).$\diamondsuit$

\vskip 0.2 true cm

\noindent{\bf Lemma 3.3.}  Suppose that  $n\geq4$.  Then for  $h\in(-\frac{1}{2},0)$, there exist polynomials $P_{2}(h)$,~$P_{1}(h)$ and $P_{0}(h)$ of $h$ with degree respectively $2n-3$, $2n-4$ and $2n-5$ such that $L(h)\bf(\sigma(h)\bf V_1(h)\bf )=0$, where
$$L(h)=P_{2}(h)\frac{{\rm d}^2}{{\rm d}h^2}+P_{1}(h)\frac{{\rm d}}{{\rm d}h}+P_{0}(h),\eqno{(3.20)}$$
and $$L(h)[{\bf \sigma(h)}{\bf V_1(h)}+{\bf \tau(h)}{\bf V_2(h)}]=R(h),$$
with $$R(h)=\frac{1}{h(2h+1)}\big(Q_{1}(h){\bf I_{1,0}}'(h)+Q_{2}(h){\bf I_{0,0}}'(h)\big)+Q_{3}(h){\bf I_{0,2}}'(h),\eqno{(3.21)}$$
and $Q_{i}(h)(i=1,2,3)$ are polynomials of $h$ with ${\rm deg}Q_{1}(h),{\rm deg}Q_{2}(h)\leq 3n-5$ and ${\rm deg}Q_{3}(h)\leq 3n-7$.

\vskip 0.2 true cm

\noindent{\bf Proof.} Denote by
 $$\Phi_1(h)={\bf \sigma(h)}{\bf V_1(h)},~\Phi_2(h)={\bf \tau(h)}{\bf V_2(h)}.\eqno{(3.22)}$$
We first assert that
$$\left\{\begin{array}{lc}
\Phi_1(h)=F_{n-1}(h){\bf I_{0,1}}''(h)+F_{n-1}(h){\bf I_{-1,1}}''(h),\\
\Phi_1'(h)=F_{n-2}(h){\bf I_{0,1}}''(h)+F_{n-2}(h){\bf I_{-1,1}}''(h),\\
\Phi_1''(h)=F_{n-3}(h){\bf I_{0,1}}''(h)+F_{n-3}(h){\bf I_{-1,1}}''(h),
\end{array}\right.\eqno{(3.23)}$$
where $F_{n-1}(h)$ denotes a polynomial of h with degree $n-1$ and etc. In fact, it follows from (3.17) that
$${\bf V_1'(h)}=(E-A_{1})^{-1}(A_{1}h+B_{1}){\bf V_1''(h)},$$
where $E$ is a $2\times2$ identity matrix. Hence,
$${\begin{split}
\Phi_1(h)&={\bf \sigma(h)}{\bf V_1(h)}={\bf\sigma(h)}(A_{1}h+B_{1}){\bf V_1'(h)}\\
&={\bf\sigma(h)}(A_{1}h+B_{1})(E-A_{1})^{-1}(A_{1}h+B_{1}){\bf V_1''(h)}\\
&:=F_k(h){\bf I_{0,1}}''(h)+F_l(h){\bf I_{-1,1}}''(h),
\end{split}}$$
where ${\bf\sigma(h)}=(\alpha_1(h),\alpha_2(h))$, ${\rm deg}\alpha_1(h)\leq n-3$ and ${\rm deg}\alpha_2(h)\leq n-3$. So we can calculate $k\leq n-1$ and $l\leq n-1$. For $\Phi_1'(h)$, we have
$${\begin{split}
\Phi_1'(h)&={\bf \sigma'(h)}{\bf V_1(h)}+{\bf \sigma(h)}{\bf V_1'(h)}\\
&=\left[{\bf \sigma'(h)}(A_{1}h+B_{1})+{\bf \sigma(h)}\right](E-A_{1})^{-1}(A_{1}h+B_{1}){\bf V_1''(h)}.
\end{split}}$$
The result for $\Phi_1''(h)$ can be proved similarly.

Next, suppose that
$$P_{2}(h)=\sum\limits_{k=0}^{2n-3}p_{2,k}h^k,~~P_{1}(h)=\sum\limits_{j=0}^{2n-4}p_{1,j}h^j,~~P_{0}(h)=\sum\limits_{l=0}^{2n-5}p_{0,l}h^l\eqno{(3.24)}$$
are polynomials of h with coefficients $p_{2,k},~p_{1,j}$ and $p_{0,l}$ to be determined such that $L(h)\Phi_1(h)=0$ for
$$0\leq k\leq 2n-3,~0\leq j\leq 2n-4,~0\leq l\leq 2n-5.\eqno{(3.25)}$$
By calculation, we have
$${\begin{split}
L(h)\Phi_1(h)&=P_{2}(h)\Phi_1''(h)+P_{1}(h)\Phi_1'(h)+P_{0}(h)\Phi_1(h)\\
&=P_{2}(h)\big[F_{n-3}(h){\bf I_{0,1}}''(h)+F_{n-3}(h){\bf I_{-1,1}}''(h)\big]\\
&\quad+P_{1}(h)\big[F_{n-2}(h){\bf I_{0,1}}''(h)+F_{n-2}(h){\bf I_{-1,1}}''(h)\big]\\
&\quad+P_{0}(h)\big[F_{n-1}(h){\bf I_{0,1}}''(h)+F_{n-1}(h){\bf I_{-1,1}}''(h)\big]\\
&=\big[P_{2}(h)F_{n-3}(h)+P_{1}(h)F_{n-2}(h)+P_{0}(h)F_{n-1}(h)\big]{\bf I_{0,1}}''(h)\\
&\quad+\big[P_{2}(h)F_{n-3}(h)+P_{1}(h)F_{n-2}(h)+P_{0}(h)F_{n-1}(h)\big]{\bf I_{-1,1}}''(h)\\
&:=X(h){\bf I_{0,1}}''(h)+Y(h){\bf I_{-1,1}}''(h)
\end{split}}$$
where $X(h)$ and $Y(h)$ are polynomials of h with degree no more than $3n-6$ and $3n-6$ respectively. Let
$$X(h)=\sum\limits_{i=0}^{3n-6}x_ih^i,~~Y(h)=\sum\limits_{j=0}^{3n-6}y_jh^j,$$
where $x_i$ and $y_j$ are expressed by $p_{2,k}$, $p_{1,j}$ and $p_{0,l}$ of (3.24) linearly, k, j and l satisfy (3.25). So $L(h)\Phi_1(h)=0$ is satisfied if we let
$$x_i=0,~~y_j=0,~~~\big(0\leq i\leq3n-6,0\leq j\leq 3n-6\big).\eqno{(3.26)}$$
System (3.26) is a homogeneous linear equation with $6n-10$ equations about $6n-9$ variables of $p_{2,k},~p_{1,j}$ and $p_{0,l}$ for $k,j,l$ satisfying (3.25). Since $(6n-9)-(6n-10)=1>0$, it follows from the theory of linear algebra that there are $p_{2,k},~p_{1,j}$ and $p_{0,l}$ such that (3.26) holds, which yields the result.

Similar to the proof of (3.23), by the second equation of (3.17) and (3.18), we can get
$${\begin{cases}
        \Phi_2(h)=F_{n-2}(h){\bf I_{1,0}}'(h)+F_{n-2}(h){\bf I_{0,0}}'(h)+F_{n-2}(h){\bf I_{0,2}}'(h),\\
        \Phi_2'(h)=F_{n-3}(h){\bf I_{1,0}}'(h)+F_{n-3}(h){\bf I_{0,0}}'(h)+F_{n-3}(h){\bf I_{0,2}}'(h),\\
        \Phi_2''(h)=\frac{1}{h(2h+1)}\left[F_{n-2}(h){\bf I_{1,0}}'(h)+F_{n-2}(h){\bf I_{0,0}}'(h)\right]+F_{n-4}(h){\bf I_{0,2}}'(h).
\end{cases}}\eqno{(3.27)}$$
Therefore, we have
$${\begin{split}
L(h)\big(\Phi_1(h)+\Phi_2(h)\big)&=L(h)\Phi_2(h)\\
&=P_{2}(h)\Phi_2''(h)+P_{1}(h)\Phi_2'(h)+P_{0}(h)\Phi_2(h)
\end{split}}\eqno{(3.28)}$$
Substituting (3.27) into (3.28), we have (3.21).$\diamondsuit$

\vskip 0.2 true cm

For $h\in(-\frac{1}{2},0)$, it is easy to calculate that
$${\bf I_{0,1}}(h)=\int_{\Gamma_h^+}x^{-4}y{\rm d}x=\iint\limits_Dx^{-4}{\rm d}x{\rm d}y\neq0.\eqno(3.29)$$
In the following discussion, we denote by $\#\{\varphi(h)=0,h\in(a,b)\}$ the number of isolated
zeros of $\varphi(h)$ on (a,b) (taking into account the multiplicity).

\vskip 0.2 true cm

\noindent{\bf Lemma 3.4.} For $h\in(-\frac{1}{2},0)$,  $\Phi_1(h)$  has at most $3n-7$ zeros (taking into account the multiplicity).

\vskip 0.2 true cm

\noindent{\bf Proof.}   Let
$$S(h)=\frac{\Phi_1(h)}{{\bf I_{0,1}}(h)}=\alpha_1(h)+\alpha_2(h)\omega(h),$$
where $\omega(h)=\frac{{\bf I_{-1,1}}(h)}{{\bf I_{0,1}}(h)}$.  It is easy to prove that $\omega(h)$ and $S(h)$ satisfy the following Riccati equation
$$h(2h+1)\omega'(h)=-\frac{2}{3}\omega^2(h)+\frac{1}{3}(4h+9)\omega(h)-\frac{1}{3}(8h+9),\eqno{(3.30)}$$
$$\alpha_2(h)h(2h+1)S'(h)=-\frac{2}{3}S(h)^2+N_{1}(h)S(h)+N_{2}(h),\eqno{(3.31)}$$
respectively, where $N_{i}(h)~(i=1,2)$ are polynomials of h with ${\rm deg}N_{1}(h)\leq n-2$ and ${\rm deg}N_{2}(h)\leq 2n-5$. In fact, using first equation of (3.17) and by direction computation, we have (3.30). From (3.30), we have
$${\begin{split}
h(2h+1)\alpha_2(h)S'(h)&=h(2h+1)\alpha_2(h)\big(\alpha_1'(h)+\alpha_2'(h)\omega(h)+\alpha_2(h)\omega'(h)\big)\\
&=h(2h+1)\alpha_2(h)\alpha_1'(h)+h(2h+1)\alpha_2(h)\alpha_2'(h)\omega(h)\\
&\quad+\alpha_2(h)^2\big(-\frac{2}{3}\omega^2(h)+\frac{1}{3}(4h+9)\omega(h)-\frac{1}{3}(8h+9)\big)\\
&=-\frac{2}{3}S(h)^2+\big(h(2h+1)\alpha_2'(h)+(\frac{4}{3}h+3)\alpha_2(h)\big)S(h)\\
&\quad+h(2h+1)\alpha_1'(h)\alpha_2(h)-h(2h+1)\alpha_1(h)\alpha_2'(h)-\frac{2}{3}\alpha_1(h)^2\\
&\quad-(\frac{4}{3}h+3)\alpha_1(h)\alpha_2(h)-(\frac{8}{3}h+3)\alpha_2(h)^2.
\end{split}}$$
Since ${\rm deg}\alpha_1(h),{\rm deg}\alpha_2(h)\leq n-3$ , we get (3.31). We can prove the desired result now. By (3.31) and Lemma 4.4 of [20], we can obtain that
$${\begin{split}
\#\{\Theta_1(h)=0,h\in(-\frac{1}{2},0)\}&=\#\{S(h)=0,h\in(-\frac{1}{2},0)\}\\
&\leq\#\{\alpha_2(h)=0,h\in(-\frac{1}{2},0)\}+\#\{N_{2}(h),h\in(-\frac{1}{2},0)\}+1\\
&\leq n-3+2n-5+1\\
&=3n-7.
\end{split}}$$
This ends the proof.$\diamondsuit$

\vskip 0.2 true cm
\noindent{\bf Lemma 3.5.} The function $R(h)$ ($h\in(-\frac{1}{2},0)$) has at most $21n-37$ zeros (taking into account the multiplicity).

\vskip 0.2 true cm

\noindent{\bf Proof.}  By (3.18) and (3.21), we can obtain
$$R^{(3n-6)}(h)=\frac{1}{[h(2h+1)]^{3n-5}}\big(\widetilde{Q}_{1}(h){\bf I_{0,0}}'(h)+\widetilde{Q}_{2}(h){\bf I_{1,0}}'(h)\big),\eqno{(3.32)}$$
where $\widetilde{Q}_{i}(h)~(i=1,2)$ are polynomials of h with degree no more than $6n-11$.
Since ${\bf I_{0,0}}(h)=\frac{1}{3}\big(\frac{1}{x_a(h)^3}-\frac{1}{x_b(h)^3}\big)$ strictly increase for $h\in(-\frac{1}{2},0)$, we have that ${\bf I_{0,0}}'(h)\neq0$. Let $\chi(h)=\frac{{\bf I_{1,0}}'(h)}{{\bf I_{0,0}}'(h)}$, by (3.18), we have
$$h(2h+1)\chi'(h)=\frac{1}{18}\big(800h+513\big)\chi(h)^2-\frac{1}{9}\big(62h+\frac{171}{2}\big)\chi(h)+4h.$$
Set
$$T(h)=\widetilde{Q}_{1}(h)+\widetilde{Q}_{2}(h)\chi(h),$$
following the line of the proof of (3.31), we can prove that
$$h(2h+1)\widetilde{Q}_{2}(h)T'(h)=\frac{1}{2}\big(800h+513\big)T(h)^2+\widetilde{N}_{1}(h)T(h)+\widetilde{N}_{2}(h),$$
where $\widetilde{N}_{1}(h)$ and $\widetilde{N}_{2}(h)$ are polynomials of h with degree no more than $6n-10$ and $12n-21$, respectively.

Hence, using Lemma 4.4 of [20] and Roll theorem, we have
$${\begin{split}
\#\{R(h)=0,h\in(-\frac{1}{2},0)\}&\leq \#\{T(h)=0,h\in(-\frac{1}{2},0)\}+3n-6\\
&\leq\#\{\widetilde{Q}_{2}(h)=0,h\in(-\frac{1}{2},0)\}+\#\{\widetilde{N}_{2}(h)=0,h\in(-\frac{1}{2},0)\}\\
&\quad+1+3n-6\\
&\leq 6n-11+12n-21+1+3n-6\\
&=21n-37,
\end{split}}$$
which ends the proof.$\diamondsuit$

\vskip 0.2 true cm

\noindent{\bf Proof of Theorem 1.1 for system (1.5)} If  $n\geq4$, then  $M(h)=\frac{1}{h^{n-3}}(\Phi_1(h)+\Phi_2(h))$,
and $L(h)(\Phi_1(h)+\Phi_2(h))=R(h)$, where $L(h)$ and $R(h)$ are given by (3.20) and (3.21), $h\in(-\frac{1}{2},0)$. It follows from
Lemma 3.4 that $\Phi_1(h)$ has at most $3n-7$ zeros on $(-\frac{1}{2},0)$. We assume that
$$P_{2}(\tilde{h}_i)=0,~\Phi_1(\bar{h}_j)=0,~~\tilde{h}_i,\bar{h}_j\in(-\frac{1}{2},0),1\leq i\leq2n-3, 1\leq j\leq 3n-7.$$
Denote $\tilde{h}_i$ and $\bar{h}_j$ as $h_m^*$, and reorder them such that $h_m^*<h_{m+1}^*$ for $m=1,\ldots,5n-10$. Let
$$\Sigma_s=(h_s^*,h_{s+1}^*),~~~s=0,1,\ldots,5n-10,$$
where $h_0^*=-\frac{1}{2}$ and $h_{5n-9}^*=0$. Then $P_{2}(h)\neq0$ and $\Phi_1(h)\neq0$ for $h\in\Sigma_s$, and the solution space of
$$L(h)=P_{2}(h)\left(\frac{{\rm d}^2}{{\rm d}h^2}+\frac{P_{1}(h)}{P_{2}(h)}\frac{{\rm d}}{{\rm d}h}+\frac{P_{0}(h)}{P_{2}(h)}\right)$$
is a Chebyshev space on $(-\frac{1}{2},0)$. By proposition 2.2, $\Phi_1(h)+\Phi_2(h)$ has at most $2+l_s$ zeros for $h\in\Sigma_s$, where $l_s$ is the number of zeros of $R(h)$ on $\Sigma_s$. Hence, by Lemma 3.5, we have
$${\begin{split}
\#\{M(h)=0,h\in(-\frac{1}{2},0)\}&=\#\{\Phi_1(h)+\Phi_2(h)=0,h\in(-\frac{1}{2},0)\}\\
&\leq\#\{R(h)=0,h\in(-\frac{1}{2},0)\}+2\cdot \emph{the number of intervals of }\Sigma_s\\
&\quad+\emph{the number of the end points of }\Sigma_s\\
&\leq 21n-37+2(5n-10+1)+5n-10\\
&=36n-65.
\end{split}}$$
When $n=1,2,3$, we can get the results of Theorem 1.1 for system (1.5) following the above line.

\vskip 0.2 true cm
\centerline {\bf \large 4 Proof of Theorem 1.1: for system (1.6)}
\vskip 0.2 true cm

In this section, we will prove Theorem 1.1 for system (1.6). For $\varepsilon=0$, the first integral of system (1.6) is
$$H(x,y)=\frac{1}{2}y^2+x-\frac{1}{3}x^3=h,~~~h\in(-\frac{2}{3},\frac{2}{3}).\eqno{(4.1)}$$
The points $(-1,0)$ and $(1,0)$ are  elementary center and saddle  corresponding to $h=-\frac{2}{3}$ and $h=\frac{2}{3}$, respectively (see Fig.\,3).

Then by (2.4), the first order Melnikov function for system (1.6) is
$$M(h)=\int_{\Gamma_h^+}g^+(x,y){\rm d}x-f^+(x,y){\rm d}y+\int_{\Gamma_h^-}g^-(x,y){\rm d}x-f^-(x,y){\rm d}y,\eqno{(4.2)}$$
where
$$\Gamma_h^+(\Gamma_h^-)=\{\left.(x,y)\right|H(x,y)=h,h\in(-\frac{2}{3},\frac{2}{3}),y>0(y<0)\}.$$

{\bf For convenience,  in this section, we denote by}
$${\bf I_{i,j}(h)=\int_{\Gamma_h^+}x^iy^j{\rm d}x}.$$

\vskip 0.2 true cm

\noindent{\bf Lemma 4.1.} Suppose that $h\in(-\frac{2}{3},\frac{2}{3})$, then we have that
$$M(h)=\left[\beta_1(h){\bf I_{0,0}}(h)+\beta_2(h){\bf I_{1,0}}(h)\right]+\left[\beta_3(h){\bf I_{0,1}}(h)+\beta_4(h){\bf I_{1,1}}(h)\right]
\eqno{(4.3)}$$
where $\beta_1(h),\beta_2(h),\beta_3(h)$ and $\beta_4(h)$ are polynomials of $h$, and
$${\rm deg}\beta_1(h)\leq \left[\frac{n}{2}\right],~{\rm deg}\beta_2(h),{\rm deg}\beta_3(h)\leq \left[\frac{n-1}{2}\right],~{\rm deg}\beta_4(h)\leq \left[\frac{n-2}{2}\right].$$

\vskip 0.2 true cm

\noindent{\bf Proof.} Similar to the proof of (3.8), we can obtain
$$M(h)=\sum\limits_{i+j=0}^n\rho_{i,j}I_{i,j}(h),\eqno{(4.4)}$$
where $\rho_{i,j}$ are arbitrary real constants.

Differentiating (4.1) with respect to $x$, we get
$$y\frac{\partial y}{\partial x}+1-x^2=0.\eqno{(4.5)}$$
Multiplying $(4.5)$ by $x^{i-2}y^j{\rm d}x$ $(i\geq2)$, integrating over $\Gamma_h^+$, we have
$$\int_{\Gamma_h^+}x^{i-2}y^{j+1}{\rm d}y+\int_{\Gamma_h^+}x^{i-2}y^j{\rm d}x-\int_{\Gamma_h^+}x^iy^j{\rm d}x=0.\eqno{(4.6)}$$
By Green's formula, we have
$$\int_{\Gamma_h^+}y^j{\rm d}y=\oint_{\Gamma_h^+\cup\overrightarrow{BA}}y^j{\rm d}y=-\iint0{\rm d}x{\rm d}y.\eqno{(4.7)}$$
Hence, from (4.6) and (4.7), it holds that
$$I_{2,j}(h)=I_{0,j}(h).\eqno{(4.8)}$$
Using Green's formula, we can obtain that
$$\int_{\Gamma_h^+}x^iy^j{\rm d}y=-\frac{i}{j+1}\int_{\Gamma_h^+}x^{i-1}y^{j+1}{\rm d}x,~~(i\geq1).\eqno{(4.9)}$$
If $i\geq3$, by (4.9) and (4.8), we have
$$I_{i,j}(h)=I_{i-2,j}(h)-\frac{i-2}{j+2}I_{i-3,j+2}(h),~~~(i\geq3).\eqno{(4.10)}$$
Similarly, multiply (4.1) by $x^iy^{j-2}{\rm d}x~~(j\geq2)$ and integrating over $\Gamma_h^+$ yields
$$I_{i,j}(h)=2hI_{i,j-2}(h)-2I_{i+1,j-2}(h)+\frac{2}{3}I_{i+3,j-2}(h),~~~(j\geq2).\eqno{(4.11)}$$
It follows from $(4.10)$ and $(4.11)$ that
$$I_{i,j}(h)=\frac{6j}{2i+3j+2}\left(hI_{i,j-2}(h)-\frac{2}{3}I_{i+1,j-2}(h)\right),~~~(j\geq2).\eqno{(4.12)}$$
Taking $j=0$ in (4.8), $(i,j)=(0,2)$ in (4.12), respectively, we have
$$I_{2,0}(h)={\bf I_{0,0}}(h),~{\bf I_{0,2}}(h)=\frac{3}{2}h{\bf I_{0,0}}(h)-{\bf I_{1,0}}(h).\eqno{(4.13)}$$
Taking $(i,j)=(3,0)$ in (4.10), $j=1$ in (4.8), and $(i,j)=(1,2),(0,3)$ in (4.12), respectively, we get
$${\begin{split}
&I_{3,0}(h)=-\frac{3}{4}h{\bf I_{0,0}}(h)+\frac{3}{2}{\bf I_{1,0}}(h),~~I_{2,1}(h)={\bf I_{0,1}}(h),\\
&I_{1,2}(h)=-\frac{4}{5}{\bf I_{0,0}}(h)+\frac{6}{5}h{\bf I_{1,0}}(h),~I_{0,3}(h)=\frac{18}{11}h{\bf I_{0,1}}(h)-\frac{12}{11}{\bf I_{1,1}}(h).
\end{split}}\eqno{(4.14)}$$

Now we prove the statement by induction on $n$. In fact, $(4.13)$ and $(4.14)$ imply that the statement holds for $n=2,3$. Now assume that it is true for $i+j\leq k-1~(k\geq3)$, then for $i+j=k$ taking $(i,j)=(0,k),(1,k-1)$ in $(4.12)$, $j=k-2$ in $(4.8)$ and $(i,j)=(3,k-3),\ldots,(k,0)$ in $(4.10)$, respectively, we have
$$
\left(\begin{matrix}
                I_{0,k}(h)\\
                 I_{1,k-1}(h)\\
                 I_{2,k-2}(h)\\
                 I_{3,k-3}(h)\\
                 \vdots\\
                  I_{k-1,1}(h)\\
                  I_{k,0}(h)
                \end{matrix}\right)\ \
=\left(\begin{matrix}
                \frac{6k}{3k+2}hI_{0,k-2}(h)-\frac{4k}{3k+2}I_{1,k-2}(h)\\
               \frac{6(k-1)}{3k+1}hI_{1,k-3}(h)-\frac{4(k-1)}{3k+1}I_{2,k-3}(h)\\
                I_{0,k-2}(h)\\
                I_{1,k-3}(h)-\frac{1}{k-1}I_{0,k-1}(h)\\
                                  \vdots\\
                  I_{k-3,1}(h)-\frac{k-3}{3}I_{k-4,3}(h)\\
                I_{k-2,0}(h)-\frac{k-2}{2}I_{k-3,2}(h)
                \end{matrix}\right).\eqno{(4.15)}
$$

By the induction hypothesis we obtain the expression $(4.3)$. Next, we estimate the degree of polynomials $\beta_i(h)~(i=1,\ldots,4)$. From (4.15), if $(i,j)=(0,k),(1,k-1)$, we have
 $${\begin{split}
 I_{i,j}(h)&=h\left[\alpha^{(k-2)}(h){\bf I_{0,0}}(h)+\beta^{(k-2)}(h){\bf I_{1,0}}(h)+\gamma^{(k-2)}(h){\bf I_{0,1}}(h)+\delta^{(k-2)}(h)I_{1,1}(h)\right]\\
 &\quad+\left[\alpha^{(k-1)}(h){\bf I_{0,0}}(h)+\beta^{(k-1)}(h){\bf I_{1,0}}(h)+\delta^{(k-1)}(h){\bf I_{0,1}}(h)+\delta^{(k-1)}(h)I_{1,1}(h)\right]\\
 &:=\alpha^{(k)}(h){\bf I_{0,0}}(h)+\beta^{(k)}(h){\bf I_{1,0}}(h)+\gamma^{(k)}(h){\bf I_{0,1}}(h)+\delta^{(k)}(h)I_{1,1}(h),
 \end{split}}$$
where $\alpha^{(k-s)}(h),~\beta^{(k-s)}(h),~\gamma^{(k-s)}(h)$ and $\delta^{(k-s)}(h)$ ($s=1,2$) are polynomials in $h$ and satisfy
$${\rm deg}\alpha^{(k-s)}(h)\leq \left[\frac{k-s}{2}\right],~{\rm deg}\beta^{(k-s)}(h),{\rm deg}\gamma^{(k-s)}(h)\leq\left[\frac{k-s-1}{2}\right],~{\rm deg}\delta^{(k-s)}(h)\leq\left[\frac{k-s-2}{2}\right],$$
$s=1,2.$ Therefore, we have
$${\rm deg}\alpha^{(k)}(h)\leq\left[\frac{k}{2}\right],~{\rm deg}\beta^{(k)}(h),{\rm deg}\gamma^{(k)}(h)\leq\left[\frac{k-1}{2}\right],~{\rm deg}\delta^{(k)}(h)\leq\left[\frac{k-2}{2}\right].$$
For other $(i,j)~(i+j=k)$ we can estimate the degree of $\beta_i(h)~(i=1,\ldots,4)$, similarly. This ends the proof.$\diamondsuit$

\vskip 0.2 true cm

\noindent{\bf Lemma 4.2.} The vector functions $\overline{{\bf V}}_1(h)=\left({\bf I_{0,0}}(h),{\bf I_{1,0}}(h)\right)^T$ and $\overline{{\bf V}}_2(h)=\left({\bf I_{0,1}}(h),{\bf I_{1,1}}(h)\right)^T$ satisfy  the following Picard-Fuchs equations
$$\overline{{\bf V}}_1(h)=(\overline{A}_1h+\overline{B}_1)\overline{{\bf V}}_1'(h),\eqno{(4.16)}$$
$$\overline{{\bf V}}_2(h)=(\overline{A}_2h+\overline{B}_2)\overline{{\bf V}}_2'(h),\eqno{(4.17)}$$
respectively, where
$$\overline{A}_1=\left(\begin{matrix}
                3&0\\
                 0&\frac{3}{2}
                \end{matrix}\right), ~~~~~
            \overline{B}_1=\left(\begin{matrix}
                0&-2\\
                -1&0
                \end{matrix}\right),
$$
$$\overline{A}_2=\left(\begin{matrix}
                \frac{6}{5}&0\\
                 0&\frac{6}{7}
                \end{matrix}\right), ~~~~~~
            \overline{B}_2=\left(\begin{matrix}
                0&-\frac{4}{5}\\
                -\frac{4}{7}&0
                \end{matrix}\right).
$$

\vskip 0.2 true cm

\noindent{\bf Proof.}
Suppose that $\Gamma_h^+$ intersects with x-axis at $(x_a(h),0)$ and $(x_b(h),0)$, then
$$I_{i,j}'(h)=\int_{x_a(h)}^{x_b(h)}jxy^{j-1}{\rm d}x+x_b(h)^iy(x_b(h),h)^j\frac{\partial{x_b(h)}}{\partial{h}}-x_a(h)^iy(x_a(h),h)^j\frac{\partial{x_a(h)}}{\partial{h}}.$$
By $x_b(h)-\frac{1}{3}x_b(h)^3=h$, we can obtain that
$$\frac{\partial{x_b(h)}}{\partial{h}}=\frac{1}{1-x_b(h)^2}.$$
Since $x_b(h)\in(-2,-1)$, we have $\frac{\partial{x_b(h)}}{\partial{h}}\neq\infty$. Similarly, $\frac{\partial{x_a(h)}}{\partial{h}}\neq\infty$. Therefore, by $y(x_a(h),h)=y(x_b(h),h)=0$, we have
$$I_{i,j}'(h)=\int_{x_a(h)}^{x_b(h)}jxy^{j-1}{\rm d}x.$$
From (4.1), we can obtain that
$$\frac{\partial y}{\partial h}=\frac{1}{y},$$
which implies
$$I_{i,j}'(h)=j\int_{\Gamma_h^+}x^iy^{j-2}{\rm d}x=jI_{i,j-2}(h).$$
Therefore, we have
$$I_{i,j}(h)=\frac{1}{j+2}I'_{i,j+2}(h).\eqno{(4.18)}$$
From $(4.18)$ and noting that ${\bf I_{0,2}}(h)=\frac{3}{2}h{\bf I_{0,0}}(h)-{\bf I_{1,0}}(h)$, we have
$${\bf I_{0,0}}(h)=3h{\bf I_{0,0}}'(h)-2{\bf I_{1,0}}'(h).$$
Similarly, we can obtain
$${\begin{split}
{\bf I_{1,0}}(h)&=-{\bf I_{0,0}'}(h)+\frac{3}{2}h{\bf I_{1,0}}'(h),\\
{\bf I_{0,1}}(h)&=\frac{6}{5}h{\bf I_{0,1}}'(h)-\frac{4}{5}{\bf I_{1,1}}'(h),\\
{\bf I_{1,1}}(h)&=-\frac{4}{7}{\bf I_{0,1}}'(h)+\frac{6}{7}h{\bf I_{1,1}}'(h).
\end{split}}$$
This ends the proof.$\diamondsuit$

\vskip 0.2 true cm

Using (4.16) and (4.17), by direct calculation, we have the following results.

\noindent{\bf Lemma 4.3.} For $h\in(-\frac{2}{3},\frac{2}{3})$, we have the following results.\\
(1) $$\overline{{\bf V}}_1(h)=(-\frac{9}{2}E_2h^2+\overline{C}_1h+\overline{D}_1)\overline{{\bf V}}_1''(h),\eqno{(4.19)}$$
$$\overline{{\bf V}}_1'(h)=(\overline{C}_2h+\overline{D}_2)\overline{{\bf V}}_1''(h),\eqno{(4.20)}$$
where $E_2$ is the identity matrix and
$$\overline{C}_1=\left(\begin{matrix}
                0&9\\
                 \frac{9}{2}&0
                \end{matrix}\right), ~~~~~
            \overline{D}_1=\left(\begin{matrix}
                -4&0\\
                0&-1
                \end{matrix}\right),
$$
$$\overline{C}_2=\left(\begin{matrix}
                -\frac{3}{2}&0\\
                 0&-3
                \end{matrix}\right), ~~~~~
            \overline{D}_2=\left(\begin{matrix}
                0&1\\
                2&0
                \end{matrix}\right).
$$
(2) $$\overline{{\bf V}}_2(h)=(\overline{C}_3h^2+\overline{D}_3)\overline{{\bf V}}_2''(h),\eqno{(4.21)}$$
$$\overline{{\bf V}}_2'(h)=(\overline{C}_4h+\overline{D}_4)\overline{{\bf V}}_2''(h),\eqno{(4.22)}$$
where
$$\overline{C}_3=\left(\begin{matrix}
                -\frac{36}{5}&0\\
                 0&\frac{36}{7}
                \end{matrix}\right), ~~~~~
            \overline{D}_3=\left(\begin{matrix}
                \frac{16}{5}&0\\
                0&-\frac{16}{7}
                \end{matrix}\right),
$$
$$\overline{C}_4=\left(\begin{matrix}
                -6&0\\
                 0&6
                \end{matrix}\right), ~~~~~
            \overline{D}_4=\left(\begin{matrix}
                0&4\\
                4&0
                \end{matrix}\right).
$$

\vskip 0.2 true cm

Set
$$\overline{{\bf \sigma}}(h)=(\beta_1(h),\beta_2(h)),~~~~\overline{{\bf \tau}}(h)=(\beta_3(h),\beta_4(h)),$$
$$M(h)=\overline{{\bf \sigma}}(h)\overline{{\bf V}}_1(h)+\overline{{\bf \tau}}(h)\overline{{\bf V}}_2(h) :=\Psi_1(h)+\Psi_2(h),$$
where $\beta_i(h)~(i=1,\ldots,4)$ and $\overline{{\bf V}}_i(h)~(i=1,2)$ are in Lemma 4.1 and Lemma 4.2, respectively. Following the proof of Lemma 3.4, by (4.21) and (4.22), we have the following results.

\vskip 0.2 true cm

\noindent{\bf Lemma 4.4.} For $h\in(-\frac{2}{3},\frac{2}{3})$, there exist polynomials polynomials $\overline{P}_2(h),~\overline{P}_1(h)$ and $\overline{P}_0(h)$ of $h$ with degree respectively $n+1,n$ and $n-1$ such that $\overline{L}(h)\Psi_2(h)=0$, where
$$\overline{L}(h)=\overline{P}_2(h)\frac{{\rm d}^2}{{\rm d}h^2}+\overline{P}_1(h)\frac{{\rm d}}{{\rm d}h}+\overline{P}_0(h).\eqno{(4.23)}$$

\vskip 0.2 true cm

Similar to Lemma 3.5, we have the following results.

\noindent{\bf Lemma 4.5.} $\overline{L}(h)\overline{M}(h)=\overline{L}(h)\left(\Psi_1(h)+\Psi_2(h)\right)=\overline{R}(h)$, where
$$\overline{R}(h)=\overline{Q}_1(h)I''_{0,0}(h)+\overline{Q}_2(h)I''_{1,0}(h),\eqno{(4.24)}$$
$\overline{Q}_1(h)$ and  $\overline{Q}_2(h)$ are polynomials of $h$ with ${\rm deg}\overline{Q}_1(h)\leq n+[\frac{n}{2}]+1,{\rm deg}\overline{Q}_2(h)\leq n+[\frac{n-1}{2}]+1$.

\vskip 0.2 true cm
Next, we estimate the number of zeros for $\Psi_2(h)$ and $\overline{R}(h)$ on $h\in(-\frac{2}{3},\frac{2}{3})$, then finish the proof of Theorem 1.1.

\vskip 0.2 true cm

\noindent{\bf Lemma 4.6.} For $h\in(-\frac{2}{3},\frac{2}{3})$,  $\Psi_2(h)$ has most $2[\frac{n}{2}]+[\frac{n-1}{2}]$ zeros (taking into account the multiplicity).

\vskip 0.2 true cm

\noindent{\bf Proof.} For $h\in(-\frac{2}{3},\frac{2}{3})$, we have
 $${\begin{split}
 {\bf I_{0,1}}(h)&=\int_{\Gamma_h^+}y{\rm d}x=\oint_{\Gamma_h^+\cup \overrightarrow{BA}}y{\rm d}x-\int_{\overrightarrow{BA}}y{\rm d}x\\
 &=\iint\limits_D{\rm d}x{\rm d}y\neq0.
 \end{split}}$$
By (4.17), we have
$$\overline{G}(h)\left(\begin{matrix}
                {\bf I_{0,1}}'(h)\\
                {\bf I_{1,1}}'(h)
                 \end{matrix}\right)
=\left(\begin{matrix}
                \frac{6}{7}h&\frac{4}{5}\\
                  \frac{4}{7}&\frac{6}{5}h
                \end{matrix}\right)
                \left(\begin{matrix}
                {\bf I_{0,1}}(h)\\
                 {\bf I_{1,1}}(h)
                \end{matrix}\right),\eqno{(4.25)}
$$
where $\overline{G}(h)=\frac{4}{35}(9h^2-4)$. Let $\chi(h)=\frac{{\bf I_{1,1}}(h)}{{\bf I_{0,1}}(h)}$, using (4.25), we can obtain that $\chi(h)$ satisfies the Riccati equation
$$\overline{G}(h)\chi'(h)=-\frac{4}{5}\chi(h)^2+\frac{12}{35}h\chi(h)+\frac{4}{7}.\eqno{(4.26)}$$
Let
$$\overline{S}(h):=\Psi_2(h)/{\bf I_{0,1}}(h)=\beta_2(h)+\beta_3(h)\chi(h),\eqno{(4.27)}$$
note that (4.26), we have
$${\begin{split}
\beta_4(h)\overline{G}(h)\overline{S}'(h)&=\beta_4(h)\overline{G}(h)\left(\beta_3'(h)+\beta_4'(h)\chi(h)+\beta_4(h)\chi'(h)\right)\\
&=\beta_4(h)G_1(h)\beta_3'(h)+\overline{G}(h)\beta_4'(h)\beta_4(h)\chi(h)+\beta_4(h)^2\left(-\frac{4}{5}\chi(h)^2+\frac{12}{35}h\chi(h)+\frac{4}{7}\right)\\
&=-\frac{4}{5}\overline{S}(h)^2+\left[\overline{G}(h)\beta_4'(h)+\frac{8}{5}\beta_3(h)+\frac{12}{35}h\beta_4(h)\right]\overline{S}(h)\\
&~~+\beta_4(h)\overline{G}(h)\beta_3'(h)-\beta_4'(h)\overline{G}(h)\beta_3(h)-\frac{4}{5}\beta_3(h)^2-\frac{12}{35}h\beta_4(h)\beta_3(h)+\frac{4}{7}\beta_4(h)^2.
\end{split}}$$
Since ${\rm deg}\beta_3(h)\leq [\frac{n-1}{2}]$, ${\rm deg}\beta_4(h)\leq [\frac{n-2}{2}]$ and ${\rm deg}\overline{G}(h)=2$, we have
$$\beta_4(h)\overline{G}(h)\overline{S}'(h)=-\frac{4}{5}\overline{S}(h)^2+\overline{N}_1(h)\overline{S}(h)+\overline{N}_2(h),\eqno{(4.28)}$$
where $\overline{N}_1(h)$ and $\overline{N}_2(h)$ are polynomials of $h$ with degree no more than $[\frac{n}{2}]$ and $[\frac{n}{2}]+[\frac{n-1}{2}]$.
Therefore, from (4.28) and Lemma 4.4 of [20], we can obtain that
$${\begin{split}
\#\{\Psi_2(h)=0,h\in(-\frac{2}{3},\frac{2}{3})\}&=\#\{\overline{S}(h)=0,h\in(-\frac{2}{3},\frac{2}{3})\}\\
&\leq \#\{\beta_4(h)=0,h\in(-\frac{2}{3},\frac{2}{3})\}+\#\{\overline{N}_2(h)=0,h\in(-\frac{2}{3},\frac{2}{3})\}+1\\
&\leq\left[\frac{n-2}{2}\right]+\left[\frac{n}{2}\right]+\left[\frac{n-1}{2}\right]+1\\
&=2\left[\frac{n}{2}\right]+\left[\frac{n-1}{2}\right].
\end{split}}$$ $\diamondsuit$

\vskip 0.2 true cm

\noindent{\bf Lemma 4.7.} For $h\in(\frac{2}{3},\frac{2}{3})$, $\overline{R}(h)$ has most $4n+[\frac{n-1}{2}]+5$ zeros (taking into account the multiplicity).

\vskip 0.2 true cm

\noindent{\bf Proof.} We first prove that ${\bf I_{0,0}}''(h)$ has one zero on $(-\frac{2}{3},\frac{2}{3})$.
Suppose that $\Gamma_h^+$ intersects with x-axis at $A(x_a(h),0)$ and $B(x_b(h),0)$ (see Fig.\,3), then we have the following:
$${\bf I_{0,0}}(h)=x_b(h)-x_a(h),\eqno{(4.29)}$$
$$-\frac{x_a(h)^3}{3}+x_a(h)=h,\eqno{(4.30)}$$
$$-\frac{x_b(h)^3}{3}+x_b(h)=h.\eqno{(4.31)}$$
By direct calculation of (4.30), we can obtain
$$\frac{{\rm d}^2x_a(h)}{{\rm d} h^2}=\frac{2x_a(h)}{(1-x_a(h)^2)^3}, ~~x_a(h)\in(-2,-1).\eqno{(4.32)}$$
Similarly, we have
$$\frac{{\rm d}^2x_b(h)}{{\rm d} h^2}=\frac{2x_b(h)}{(1-x_b(h)^2)^3}, ~~x_b(h)\in(-1,1).\eqno{(4.33)}$$
From (4.30), (4.31) and noting that $x_a(h)\neq x_b(h)$, we have
$$x_a(h)^2+x_b(h)x_a(h)+x_b(h)^2-3=0,$$
that is,
$$x_a(h)=-\frac{1}{2}x_b(h)-\frac{\sqrt{3}}{2}\sqrt{4-x_b(h)^2}.\eqno{(4.34)}$$
By (4.29),(4.32),(4.33) and (4.34), it holds that
$${\begin{split}
{\bf I_{0,0}}''(h)&=x''_b(h)-x''_a(h)\\
&=\frac{2x_b(h)}{(1-x_b(h)^2)^3}-\frac{2x_a(h)}{(1-x_a(h)^2)^3}\\
&=\frac{2g(x_b(h))}{(1-x_b(h)^2)^3(1-x_a(h)^2)^3},
\end{split}}$$
where
$${\begin{split}
g(x_b(h))&=x_b(h)\left(-2+\frac{1}{2}x_b(h)^2-\frac{\sqrt{3}}{2}x_b(h)\sqrt{4-x_b(h)^2}\right)^3\\
&\quad+\left(\frac{1}{2}x_b(h)+\frac{\sqrt{3}}{2}\sqrt{4-x_b(h)^2}\right)(1-x_b(h)^2)^3,~~x_b(h)\in(-1,1).
\end{split}}$$

From the curve of $\frac{g(x_b(h))}{(1+x_b(h))^3}$ for $x_b(h)\in(-1,1)$ (see Fig.\,4), we know ${\bf I_{0,0}}''(h)$ has one zero for $h\in(-\frac{2}{3},\frac{2}{3})$.

Suppose that ${\bf I_{0,0}}''(h_0)=0$, for $h\in(-\frac{2}{3},\frac{2}{3})\backslash h_0$,
let $\omega(h)=\frac{{\bf I_{1,0}}''(h)}{{\bf I_{0,0}}''(h)}$ and
$$\widetilde{R}(h):=\overline{R}(h)/{\bf I_{0,0}}''(h)=\overline{Q}_1(h)+\overline{Q}_2(h)\omega(h).\eqno{(4.35)}$$
By (4.20), we have
$$\widetilde{G}(h)\left(\begin{matrix}
                {\bf I_{0,0}}'''(h)\\
                {\bf I_{1,0}}'''(h)
                 \end{matrix}\right)
=\left(\begin{matrix}
                -\frac{3}{4}h&-\frac{2}{5}\\
                  -\frac{1}{2}&-\frac{3}{5}h
                \end{matrix}\right)
                \left(\begin{matrix}
                {\bf I_{0,0}}''(h)\\
                {\bf I_{1,0}}''(h)
                \end{matrix}\right),\eqno{(4.36)}
$$
where $\widetilde{G}(h)=\frac{1}{5}(\frac{9}{4}h^2-1)$. Therefore, $\omega(h)$ satisfies the Riccati equation
$$\widetilde{G}(h)\omega'(h)=\frac{2}{5}\omega(h)^2+\frac{3}{20}h\omega(h)-\frac{1}{2}.\eqno{(4.37)}$$
Then, by (4.37),  $\widetilde{R}(h)$ satisfies that
$$\overline{Q}_2(h)\widetilde{G}(h)\widetilde{R}'(h)=\frac{2}{5}\widetilde{R}(h)^2+\widetilde{T}_1(h)\widetilde{R}(h)+\widetilde{T}_2(h),$$
where $\widetilde{T}_1(h)$, $\widetilde{T}_2(h)$ are polynomials of $h$ with degree no more than $n+[\frac{n+1}{2}]+1$ and $3n+2$. By Lemma 4.4 of [20], we have
$${\begin{split}
\#\{\overline{R}(h)=0,h\in(-\frac{2}{3},\frac{2}{3})\backslash h_0\}&=\#\{\widetilde{R}(h)=0,h\in(-\frac{2}{3},\frac{2}{3})\backslash h_0\}\\
&\leq\#\{\overline{Q}_2(h)=0,h\in(-\frac{2}{3},\frac{2}{3})\backslash h_0\}+\#\{\overline{T}_2(h)=0,h\in(-\frac{2}{3},\frac{2}{3})\backslash h_0\}+1\\
&\leq n+\left[\frac{n-1}{2}\right]+1+3n+2+1\\
&=4n+\left[\frac{n-1}{2}\right]+4.
\end{split}}$$
Therefore, $\#\{\overline{R}(h)=0,h\in(-\frac{2}{3},\frac{2}{3})\}\leq4n+[\frac{n-1}{2}]+5.$ $\diamondsuit$

\vskip 0.2 true cm

\noindent{\bf Proof of Theorem 1.1 for system (1.6).} For $h\in (-\frac{2}{3},\frac{2}{3})$, $M(h)=\Psi_1(h)+\Psi_2(h)$, and $\overline{L}(h)M(h)=\overline{R}(h)$, where $\overline{L}(h)$ and $\overline{R}(h)$ are given by (4.23) and (4.24).
It follows from Lemma 4.4 that $\Psi_2(h)$ has at most $2[\frac{n}{2}]+[\frac{n-1}{2}]$ zeros on $(-\frac{2}{3},\frac{2}{3})$. We assume that
$$\overline{P}_2(\tilde{h}_i)=0,~~~\Psi_2(\bar{h}_j)=0,~~~~\tilde{h}_i,\bar{h}_j\in(-\frac{2}{3},\frac{2}{3}),1\leq i\leq n+1,1\leq j\leq 2\left[\frac{n}{2}\right]+\left[\frac{n-1}{2}\right].$$
Denote $\tilde{h}_i$ and $\bar{h}_j$ as $h_m^*$, and reorder them such that $h_m^*<h_{m+1}^*$ for $m=1,\ldots,n+2[\frac{n}{2}]+[\frac{n-1}{2}]$. Let
$$\Sigma_s=(h_s^*,h_{s+1}^*),~~~s=0,1,\ldots,n+2\left[\frac{n}{2}\right]+\left[\frac{n-1}{2}\right]+1,$$
where $h_0^*=-\frac{2}{3}$ and $h_{n+2[\frac{n}{2}]+[\frac{n-1}{2}]+2}^*=\frac{2}{3}$. Then $\overline{P}_2(h)\neq0$ and $\Psi_2(h)\neq0$ for $h\in\Sigma_s$, and the solution space of
$$\overline{L}(h)=\overline{P}_2(h)\left(\frac{{\rm d}^2}{{\rm d}h^2}+\frac{\overline{P}_1(h)}{\overline{P}_2(h)}\frac{{\rm d}}{{\rm d}h}+\frac{\overline{P}_0(h)}{\overline{P}_2(h)}\right)$$
is a Chebyshev space on $\Sigma_s$. By proposition 2.2, $M(h)$ has at most $2+l_s$ zeros for $h\in\Sigma_s$, where $l_s$ is the number of zeros of $\overline{R}(h)$ on $\Sigma_s$. Hence, by Lemma 4.7, we have
$${\begin{split}
\#\{M(h)=0,h\in(-\frac{2}{3},\frac{2}{3})\}&\leq\#\{\overline{R}(h)=0,h\in(-\frac{2}{3},\frac{2}{3})\}+2\cdot \mbox{\it the number of intervals of } \Sigma_s\\
&\quad+\mbox{\it the number of the end points of } \Sigma_s\\
&\leq4n+\left[\frac{n-1}{2}\right]+5+2\left(n+2\left[\frac{n}{2}\right]+\left[\frac{n-1}{2}\right]+2\right)\\
&~~+n+2\left[\frac{n}{2}\right]+\left[\frac{n-1}{2}\right]+1\\
&\leq12n+6.
\end{split}}$$

\end{document}